%% file: TJACk6.tex
\documentclass[12pt, leqno]{amsart}

\usepackage{hyperref}

\usepackage{epsfig}
\usepackage[all]{xy}
\usepackage{amsopn,amsthm,amssymb,amsmath,amssymb,amsfonts}

\theoremstyle{plain}
\newtheorem{Thm}{Theorem}[section]
\newtheorem{Prop}[Thm]{Proposition}
\newtheorem{Lem}[Thm]{Lemma}
\newtheorem{Cor}[Thm]{Corollary}

\theoremstyle{definition}
\newtheorem*{Def}{Definition}

\theoremstyle{remark}

\newcommand{\arw}[1]{\ar@{-}[#1]\ar@{>}[#1]}
\newcommand{\arwstrut}[1]{\ar@{}[#1]}
\newcommand{\aru}[2]{\arw{#1}^-{\displaystyle{#2}}}
\newcommand{\ard}[2]{\arw{#1}_-{\displaystyle{#2}}}

\newcommand{\arsub}{\ar@{}[r]|{\subset}}
\newcommand{\arcap}{\ar@{}[d]|{\bigcap}}
\newcommand{\arcup}{\ar@{}[u]|{\bigcup}}

\newcommand{\trianglerd}[6]{\xymatrix{{#1} \aru{r}{#2} \ard{dr}{#3}& {#4} \aru{d}{#5}\\& {#6}}}

\newcommand{\squarerdcomm}[8]{\xymatrix{{#1} \aru{r}{#2} \ard{d}{#3}& {#4} \aru{d}{#5}\\ {#6} \ard{r}{#7}& {#8} \ar@{}[ul]|{\circlearrowleft}}}

\newcommand{\zz}{\mathbb{Z}}
\newcommand{\rr}{\mathbb{R}}
\newcommand{\qq}{\mathbb{Q}}
\newcommand{\e}{\varepsilon}
\newcommand{\ph}{\varphi}
\newcommand{\oo}{\mathcal{O}}
\newcommand{\nnn}{\mathcal{N}}
\newcommand{\ccx}{\mathcal{X}}
\newcommand{\rmap}{\longrightarrow}

\newcommand{\dvec}[1]{\displaystyle\overrightarrow{#1}}
\newcommand{\ddvec}[1]{\overrightarrow{#1}}
\newcommand{\xxangle}{\theta}
\newcommand{\xxTN}{\mathcal{P}}

\newcommand{\xxDiv}{\operatorname{Div}}
\newcommand{\xxJac}{\operatorname{Jac}}
\newcommand{\xxT}{\operatorname{T}}
\newcommand{\xxInt}{\operatorname{Int}}
\newcommand{\xxwt}{\operatorname{wt}}
\newcommand{\xxNewt}{\operatorname{Newt}}
\newcommand{\xxConv}{\operatorname{Conv}}
\newcommand{\xxVer}{\operatorname{Ver}}
\newcommand{\xxEdge}{\operatorname{Edge}}
\newcommand{\xxmoment}{\operatorname{moment}}
\newcommand{\xxMult}{\operatorname{Mult}}
\newcommand{\xxarea}{\operatorname{area}}
\newcommand{\xxlength}{\operatorname{length}}
\newcommand{\xxKer}{\operatorname{Ker}}
\newcommand{\xxPic}{\operatorname{Pic}}

\title {Tropical Jacobians in $\mathbb{R}^2$}
\author {Shuhei Yoshitomi}
\date{}
\begin {document}

\begin {abstract}

By tropical Abel-Jacobi theorem, the Jacobian of a tropical curve is isomorphic to the Picard group. A tropical curve in $\mathbb{R}^2$ corresponds to an immersion from a tropical curve to $\mathbb{R}^2$. In this paper, we show that any principal divisor on a tropical curve is the restriction of a principal divisor on the ambient plane $\mathbb{R}^2$.

\end {abstract}

\maketitle

\section {Introduction}

The main result of this paper is comparison of the Picard group of a tropical curve $X$ embedded to $\mathbb{R}^2$ and a potentially larger group obtained by quotienting the divisor group by the equivalence generated by rational functions that extend to the ambient plane. We show that these two groups are equal. In other words, any principal divisor on $X$ is induced from tropical curves in $\mathbb{R}^2$.
\par
Let $C$ be a tropical curve in $\rr^2$. Let $X \rightarrow \rr^2$ be the corresponding immersion from a tropical curve $X$ to the affine space $\rr^2$ with no crossing points. (So, if $C$ is reduced, the immersion $X \rightarrow C$ is bijective.) The Jacobian $\xxJac( X)$ is defined in \cite{Mik2}.

\begin {Thm}[Tropical Abel-Jacobi, \cite{Mik2}] \label {Thm:Z1}
The Abel's map $\mu \colon \xxDiv^0( X) \rightarrow \xxJac( X)$ factors through the Picard group $\xxPic^0( X)$:
\[
\trianglerd{\xxDiv^0( X)}{}{ \mu}{ \xxPic^0( X)}{ \phi}{ \xxJac( X)}
\]
and $\phi$ is a bijection.
\end {Thm}

\begin {Def}
A tropical curve $C$ in $\rr^2$ is {\em reduced} if every edge is of weight $1$. $C$ is {\em smooth} if every vertex is $3$-varent and of multiplicity $1$. (Hence any smooth tropical curve is reduced.)
\end {Def}

\begin {Def}
The {\em genus} of a reduced tropical curve $C$ in $\rr^2$ is the first Betti number $b_1( C)$.
\end {Def}

\begin {figure}[ht]
\begin {center} \input {./picture/bhA3_b.tex} \end {center}
\caption {Tropical elliptic curve}
\label {Fig:A3}
\end {figure}
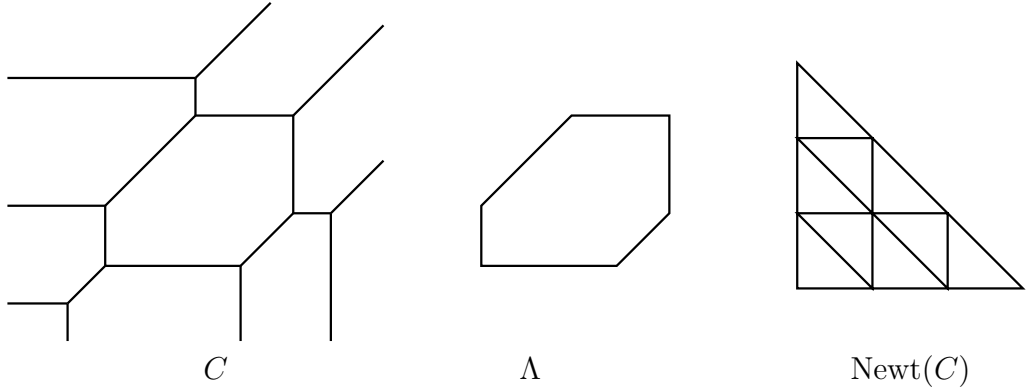

For example, a tropical elliptic curve $C$ is a smooth tropical curve of genus $1$ (Figure \ref{Fig:A3}).

\begin {Def}
A divisor $D$ of $C$ is {\em principal} if there are tropical curves $L, L'$ in $\rr^2$ such that
\[
\Delta(L)= \Delta(L'), \]\[
D= C \cdot L- C \cdot L',
\]
where $\Delta$ denotes the Newton polygon and $C \cdot L$ denotes the stable intersection. Divisors $D, D'$ of $C$ are {\em linearly equivalent} ($D \sim D'$) if $D- D'$ is principal.
\end {Def}

\begin {Thm} \label {Thm:Z2}
Let $C$ be a reduced tropical curve in $\rr^2$. Let $X \rightarrow \rr^2$ be the corresponding immersion of a tropical curve. Let $D$ be a divisor on $C$. Let $\widetilde{ D}$ be the corresponding divisor on $X$. Then $D$ is principal if and only of $\widetilde{ D}$ is principal.
\end {Thm}

Let $\xxT( C)$ be the quotient group
\[
\xxT( C)= \xxDiv^0(C) / \sim.
\]
Theorem \ref{Thm:Z2} is equivalent to the following statement.

\begin {Thm} \label {Thm:E1}
Let $C$ be a reduced tropical curve of genus $g$ in $\rr^2$. Then $\xxT( C)$ is the $g$-dimensional real torus, and $\xxT( C) \cong \xxPic^0( X)$.
\end {Thm}

Also we show the following statement in the case of genus $1$.

\begin {Thm} \label {Prop:A1}
Let $C$ be a reduced tropical curve of genus $1$ in $\rr^2$. Let $\Lambda$ be the cycle in $C$. Let $\oo$ be any point of $\Lambda$. Then the map
\begin {eqnarray*}
\ph \colon \Lambda & \rmap & \xxT( C) \\
P & \mapsto & P- \oo
\end {eqnarray*}
is bijective.
\end {Thm}

\vspace{1em}

Acknowledgements. The author thanks Professor Yujiro Kawamata for helpful advice.  He thanks the referee for his comments, which led to improvements of the paper.

\section {Preliminaries} \label {Sec:B1}

The affine space $\rr^2$ is equipped with interior product
\[
(u_1, u_2) \cdot (v_1, v_2)= u_1 v_1+ u_2 v_2
\]
and exterior product
\[
(u_1, u_2) \times (v_1, v_2)= u_1 v_2- v_1 u_2.
\]
A {\em primitive vector} in $\rr^2$ is an integral vector $u= (u_1, u_2)$ such that $u_1, u_2$ are relatively prime. Let $v \in \zz^2$ be an integral vector. Then there are a primitive vector $u$ and a number $m \in \zz_{ > 0}$ such that $v= m u$. The number $m$ is called the {\em lattice length} of $v$. We write $m= | v|$. \par
Let $C$ be a $1$-dimensional weighted simplicial complex of rational slopes in $\rr^2$ (that is a subset of $\rr^2$ written as a union of line segments equipped with weights of edges). Each finite edge $E \subset C$ has two vertices. Let $V \in E$ be one of the vertices. Let $\xxwt( E) \in \zz_{ > 0}$ be the weight of $E$. The weighted primitive vector of $E$ starting at $V$ is defined to be the integral vector $u_E$ for the direction from $V$ to $E$ such that $| u_E|= \xxwt( E)$. \par
A {\em ray} $E$ of $C$ is an infinite edge. Any ray $E$ has only one vertex.

\begin {Def}
A vertex $V$ of $C$ {\em satisfies the balancing condition} if the sum of all weighted primitive vectors starting at $V$ equals $0$:
\[
\sum_{E \ni V}^{}u_E = 0.
\]
A {\em tropical curve} in $\rr^2$ is a $1$-dimensional weighted simplicial complex of rational slopes such that each vertex satisfies the balancing condition.
\end {Def}

The union $C_1 \cup C_2$ of two tropical curves is a tropical curve. Indeed if $V \in C_1 \cap C_2$ is an intersection point of two edges, $V$ is considered as a $4$-valent vertex of $C_1 \cup C_2$. \par
Let $C$ be a tropical curve in $\rr^2$. Let $U_1, \ldots , U_r$ be all connected components of $\rr^2 \setminus C$. Let $N$ be a $1$-dimensional simplicial complex with vertex set $\xxVer(N)= \{w_1, \ldots, w_r\}$, $w_i \in \zz^2$.

\begin {Def}
$N$ is a {\em Newton complex} of $C$ if it satisfies the following conditions for any $i \not= j$. \\
i) If $\overline{U_i} \cap \overline{U_j}= \emptyset$, then $[w_i, w_j] \notin \xxEdge(N)$. \\
ii) If $\overline{U_i} \cap \overline{U_j}= E$ for some $E \in \xxEdge(C)$, then $[w_i, w_j] \in \xxEdge(N)$, and $w_j- w_i$ has lattice length $\xxwt(E)$, direction orthogonal to $E$ from $U_i$ to $U_j$.
\end {Def}

\begin {Prop} \label {Thm:A1}
Let $C, C_1, C_2$ be tropical curves in $\rr^2$. \\
1) The Newton complex $\xxNewt(C)$ of $C$ exists uniquely up to parallel translation. The convex hull
\[
\Delta(C)= \xxConv( \xxNewt(C))
\]
is called the {\em Newton polygon} of $C$. \\
2) $\Delta(C_1 \cup C_2)$ equals the Minkowski sum $\Delta(C_1)+ \Delta(C_2)$.
\end {Prop}

\begin {figure}[ht]
\begin {center} \input {./picture/bhA2.tex} \end {center}
\caption {Newton complex}
\label {Fig:A2}
\end {figure}
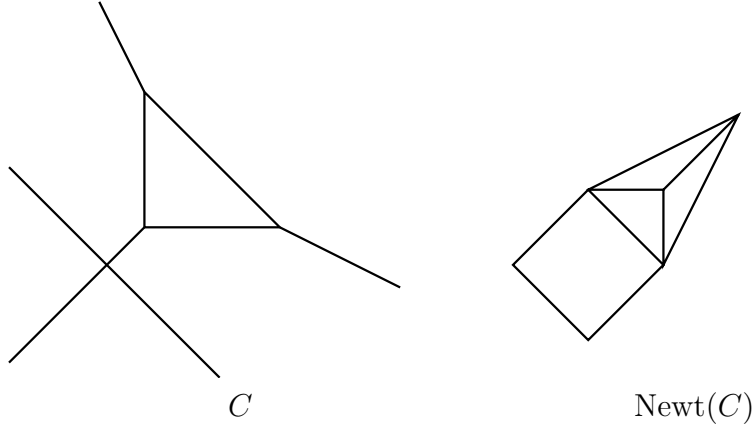

\begin {proof}
1) Let $w_1= (0, 0)$, and suppose $w_1, \ldots, w_{k-1}$ are constructed. Rearranging $U_k, \ldots, U_r$, we may assume $U_i \cap U_k= E_k$ for some $i< k$ and some $E_k \in \xxEdge(C)$. Let $u_k$ be the primitive vector of direction orthogonal to $E_k$ from $U_i$ to $U_k$. Let $w_k$ be the vector defined as follows.
\begin {equation}
w_k- w_i= \xxwt(E_k) u_k. \label {eq:T-A1-1}
\end {equation}
If $U_1, \ldots, U_k$ are adjacent at a common vertex $V$, we have the balancing condition
\[
\sum_{i=1}^{k} \xxwt(E_i) u_i= 0,
\]
where $E_1$ is the boundary of $U_k$ and $U_1$, and $u_1$ is the primitive vector of direction orthogonal to $E_1$ from $U_k$ to $U_1$. So we have
\[
w_1- w_k= \xxwt(E_1) u_1.
\]
So this construction does not depend on the choice of $U_k$.
\par
2) The set of the rays of $C_1 \cup C_2$ corresponds to the set of the rays of $C_1$ and $C_2$. So the set of the edges of $\xxNewt(C_1 \cup C_2)$ in the boundary of $\Delta(C_1 \cup C_2)$ corresponds to the set of the edges of $\xxNewt(C_1)$ and $\xxNewt(C_2)$ in the boundary of $\Delta(C_1)$ and $\Delta(C_2)$.
\end {proof}

$\xxNewt(C)$ can be considered as a dual object of $C$, with correspondence from $U_i$ to $w_i$ (Figure \ref{Fig:A2}). A vertex $V \in C$ corresponds to a polygon $T_V \subset \Delta(C)$ as follows. \\
i) $U_{i_1}, \ldots, U_{i_k}$ are adjacent at $V$, \\
ii) $T_V= \xxConv \{w_{i_1}, \ldots, w_{i_k} \}$.

\begin {Prop}[Global balancing condition] \label {Thm:A2}
Let $C$ be a tropical curve in $\rr^2$. Let $\Lambda$ be a simple closed curve in $\rr^2$ intersecting edges of $C$, say $E_1, \ldots, E_N$, transversely. Then
\[
\sum_{i=1}^{N} u_{E_i}= 0,
\]
where $u_{E_i}$ is the weighted primitive vector of $E_i$ starting at the vertex inside $\Lambda$.
\end {Prop}

\begin {proof}
For each vertex $V_j \in C$ inside $\Lambda$, the balancing condition
\[
\sum_{k} u_{jk}= 0
\]
holds. Thus
\[
\sum_{j,k} u_{jk} = 0.
\]
On the left side, two weighted primitive vectors of the same edge inside $\Lambda$ are canceled. Thus we have the required equation.
\end {proof}

A {\em tangent vector} $(v, P)$ in $\rr^2$ is a vector $v \in \rr^2$ with a starting point $P \in \rr^2$. We fix a point $P_0 \in \rr^2$. The {\em moment} of $(v, P)$ is the exterior product
\[
\xxmoment(v, P)= \ddvec{P_0 P} \times v.
\]

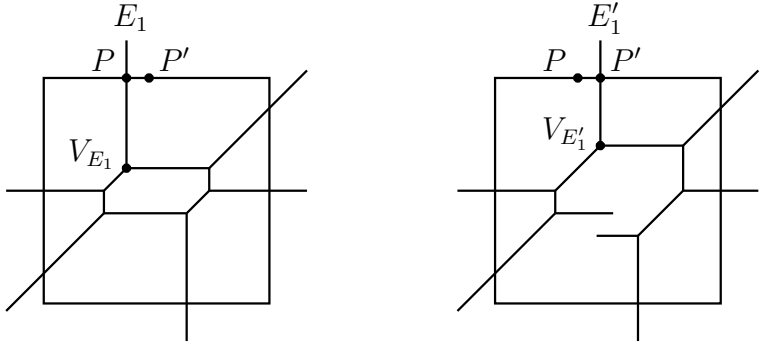
\begin {figure}[ht]
\begin {center} \input {./picture/bhB1.tex} \end {center}
\caption {Moment condition inside $\Lambda$. Note that $\xxmoment(u_{E_1}, V_{E_1}) \not= \xxmoment(u_{E'_1}, V_{E'_1})$.}
\label {Fig:B1}
\end {figure}

\begin {Prop}[Moment condition] \label {Thm:A3}
Under the assumption of Proposition \ref{Thm:A2},
\[
\sum_{i=1}^{N} \xxmoment(u_{E_i}, V_{E_i})= 0,
\]
where $(u_{E_i}, V_{E_i})$ is the weighted primitive tangent vector of $E_i$ starting at the vertex inside $\Lambda$. (See Figure \ref{Fig:B1}.)
\end {Prop}

\begin {proof}
For each vertex $V_j \in C$ inside $\Lambda$, the balancing condition
\[
\sum_{k} u_{jk}= 0
\]
holds. Thus
\[
\sum_{j,k} \ddvec{P_0 V_j} \times u_{jk}= 0.
\]
On the left side, two weighted primitive vectors of the same edge inside $\Lambda$ are canceled as follows.
\begin {eqnarray*}
\ddvec{P_0 V_j} \times u_{jk}+ \ddvec{P_0 V_{j'}} \times u_{j'k'} & = & \ddvec{P_0 V_j} \times u_{jk}+ \ddvec{P_0 V_{j'}} \times (- u_{jk}) \\
& = & \ddvec{V_{j'} V_j} \times u_{jk} \\
& = & 0.
\end {eqnarray*}
Thus we have the required equation.
\end {proof}

\begin {Def}
The {\em multiplicity} of a vertex $V \in C$ is
\[
\xxMult(V;C)= 2 \cdot \xxarea(T_V).
\]
The {\em intersection multiplicity} of $V \in C_1 \cap C_2$ is
\[
\mu_V= \frac{1}{2} \left( \xxMult(V;C_1 \cup C_2)- \xxMult(V;C_1)- \xxMult(V;C_2) \right).
\]
(If $V$ is not a vertex of $C$, we put $\xxMult(V;C)= 0$.) \\
We have a divisor
\[
C_1 \cdot C_2= \sum_{V \in C_1 \cap C_2} \mu_V V.
\]
This is called the {\em stable intersection} of $C_1$ and $C_2$.
\end {Def}

By \cite{First Steps}, the stable intersection is characterized as the limit of the transversal intersection. If $V \in C_1 \cap C_2$ is a transversal intersection point, we have
\[
\mu_V= |u_E \times u_F|,
\]
where $E \subset C_1$ and $F \subset C_2$ are edges passing through $V$.

A tropical Bezout's theorem is proved in \cite{First Steps}. Now we show another proof.

\begin {Thm}[Tropical Bezout] \label {Thm:A4}
Let $C_1, C_2$ be tropical curves in $\rr^2$. Then the following formula holds.
\[
\deg(C_1 \cdot C_2)= \xxarea( \Delta(C_1)+ \Delta(C_2))- \xxarea( \Delta(C_1))- \xxarea( \Delta(C_2)),
\]
where $\Delta(C_1)+ \Delta(C_2)$ is the Minkowski sum.
\end {Thm}

\begin {proof}
The statement follows from Proposition \ref{Thm:A1}.
\end {proof}

For example, let $c, d \geq 1$, and suppose $\Delta(C_1)= \xxConv \{(0,0), (c,0), (0,c) \}$, $\Delta(C_2)= \xxConv \{(0,0), (d,0), (0,d) \}$ ($C_2$ is said to be {\em projective} of degree $d$). Then $\Delta(C_1)+ \Delta(C_2)= \xxConv \{(0,0), (c+d,0), (0,c+d) \}$, and $\deg(C_1 \cdot C_2)= \frac{1}{2} (c+ d)^2- \frac{1}{2} c^2- \frac{1}{2} d^2= cd$. \par

\section {Reduced tropical curves} \label{Sec:F2}

Let $C$ be a reduced tropical curve.

\begin {Lem} \label {Prop:A2}
Let $u \in \zz^2$ be a primitive vector. Then for given $\e>0$, there is a primitive vector $v \in \zz^2$ such that
\[
u \times v= 1, \quad |\xxangle(u)- \xxangle(v)|< \e,
\]
where $\xxangle(u)$ denotes the angle of $u$.
\end {Lem}

This lemma is easy.

\begin {Lem} \label {Prop:A3}
Let $E$ be an edge of $C$, and let $P, P', Q, Q'$ be points of E such that $\dvec{PP'}$ $= \dvec{QQ'}$. Then $P'- P \sim Q'- Q$.
\end {Lem}

\begin {proof}
(Figure \ref{Fig:A5}) We may assume that $P', Q$ lie on the interval $[P, Q']$, and that $Q'$ lies in the interior $\xxInt(E)$. Let $v_1, v_2, v_3$ be primitive vectors such that
\begin {equation}
| u_E \times v_i |= 1 \quad (i=1, 2, 3), \label {eq:A3-1}
\end {equation}
\[
\xxangle(u_E)- \e< \xxangle(v_1)< \xxangle(u_E)< \xxangle(v_2)< \xxangle(u_E)+ \e< \xxangle(v_3).
\]
Then we have a triangle, with vertices $P$ and $R_1, R_2 \in \rr^2$, such that
\begin {eqnarray*}
&& \ddvec{PR_i} \mbox{ has direction } v_i \quad (i=1, 2), \\
&& \ddvec{R_1R_2} \mbox{ has direction } v_3, \\
&& Q' \in [R_1, R_2].
\end {eqnarray*}
Take $\e>0$ small enough so that this triangle is disjoint from $C \setminus E$. \par
    Now we take two tropical curves $L, M$ as follows. $L$ consists of one vertex $P$ and three rays $L_0, L_1, L_2 \subset \rr^2$, with $R_1 \in L_1$, $R_2 \in L_2$. $L_0$ is parallel to $E$. $\xxwt(L_1)= \xxwt(L_2)= 1$. (The balancing condition at $P$ follows from equation (\ref{eq:A3-1}).) $M$ consists of one finite edge $M_0 \subset \rr^2$, four rays $M_1, M_2, M_3, M_4 \subset \rr^2$, and two vertices $V_1, V_2 \in \rr^2$. $M_0$ is parallel to $[R_1, R_2]$, and passes through $Q$. $V_1 \in [P, R_1]$, $V_2 \in [P, R_2]$. For $i=1, 2$, $M_i$ has vertex $V_i$ and passes through $R_i$. $M_3, M_4$ are parallel to $E$. $\xxwt(M_0)= \xxwt(M_1)= \xxwt(M_2)= 1$. \par
    Move $L$ by $\dvec{PP'}$, and denote it by $L'$. Move $M$ by $\dvec{QQ'}$, and denote it by $M'$. Note that the intersection of $L'$ and $C$ at $P'$ has multiplicity $1$. Now we have
\[
(C \cdot L- P)- (C \cdot L'- P')= (C \cdot M- Q)- (C \cdot M'- Q').
\]
Thus $P'-P$ is linearly equivalent to $Q'-Q$.
\end {proof}

\begin {figure}[ht]
\begin {center} \input {./picture/bhA5.tex} \end {center}
\caption {}
\label {Fig:A5}
\end {figure}

\begin {Cor} \label {Prop:A4}
Let $E$ be any edge, and suppose that all interior points of $E$ are linearly equivalent. Then all points of $E$ are linearly equivalent.
\end {Cor}

\begin {Lem} \label {Prop:A5}
Let $E$ be a ray of $C$. Then all points of $E$ are linearly equivalent.
\end {Lem}

\begin {proof}
(Figure \ref{Fig:A6}, left) Let $P, Q \in \xxInt(E)$. Take a primitive vector $v$ so that
\[
| u_E \times v |= 1, \quad | \xxangle(u_E)- \xxangle(v) |< \e.
\]
There is a parallelogram $R_1 R_2 R_3 R_4$ such that
\begin {eqnarray*}
&& \ddvec{R_3 R_1}, \ddvec{R_4 R_2} \mbox{ have direction } v, \\
&& \ddvec{R_2 R_1}, \ddvec{R_4 R_3} \mbox{ have direction } w:= u_E- v, \\
&& P \in [R_1, R_3],~ Q \in [R_2, R_4].
\end {eqnarray*}
Take $\e>0$ small enough so that this parallelogram is disjoint from $C \setminus E$. \par
Let $M_1$ be a tropical curve, consisting of one vertex $R_1$ and three rays $L_0, L_1, L_2 \subset \rr^2$, such that $L_0$ has direction $u_E$, $R_3 \in L_1$, $R_2 \in L_2$, $\xxwt(L_0)= \xxwt(L_1)= 1$. Move $M_1$ by $\dvec{R_1 R_i}$, and denote it by $M_i (i=2, 3, 4)$. Then we have
\[
C \cdot M_1+ C \cdot M_4- Q= C \cdot M_2+ C \cdot M_3- P.
\]
Thus $P \sim Q$.
\end {proof}

\begin {figure}[ht]
\begin {center} \input {./picture/bhA6.tex} \end {center}
\caption {}
\label {Fig:A6}
\end {figure}

\begin {Def}
An edge $E \subset C$ is {\em trivial} if $C \setminus \xxInt(E)$ is disconnected.
\end {Def}

\begin {Lem} \label {Prop:A6}
Let $E$ be a trivial edge of $C$. Then all points of E are linearly equivalent.
\end {Lem}

\begin {proof}
(Figure \ref{Fig:A6}, right) Let $P \in \xxInt(E)$. Since $C \setminus \xxInt(E)$ is disconnected, $E$ is a boundary of two unbounded convex open sets $U_1, U_2$ of $\rr^2 \setminus C$. For $i=1, 2$, let $\theta_i$ be the angle of any unbounded direction of $U_i$. For $\delta, \e>0$, let
\[
W_i= \left\{
\begin {array}[c]{rcl}
P+ a+ w & | & a, w \in \rr^2,~ |a|< \delta,~ P+ a \in E, \\
&& | \xxangle(w)- \theta_i|< \e
\end {array}
\right\}, \]\[
W= W_1 \cup W_2.
\]
Take $\delta, \e>0$ small enough so that $W$ intersects $C \setminus E$ only at points of rays of $C$. \par
Let $v$ be a primitive vector such that $u_E \times v= 1$. Take $w, w' \in \zz^2$ so that
\[
w- w'= v, \]\[
| \xxangle(w)- \theta_1 |< \e, \]\[
| \xxangle(w')- \theta_1 |< \e.
\]
Let $L$ be a tropical curve consisting of three rays $L_0, L_1, L_2 \subset \rr^2$ such that
\[
L_1 \cup L_2 \subset W, \]\[
u_{L_1}= w, \]\[
| \xxangle(L_2)- \theta_2 |< \e, \]\[
L_1 \cap E= P.
\]
Moving $L$ on the direction of $L_0$, the intersection point $P$ changes to other point $Q \in E$, but all other intersection points of $L$ and $C$ are stable except for points of rays of $C$. We have
\[
(u_E \times w) P \sim (u_E \times w) Q.
\]
Similarly,
\[
(u_E \times w') P \sim (u_E \times w') Q.
\]
Thus $P \sim Q$.
\end {proof}

\section {Parameter space of tropical plane curves} \label {Sec:B3}

Let $L$ be a tropical curve with Newton complex $\nnn$. We fix a vertex $V_0= (b_1, b_2)$ of $L$. Let $E_1, \ldots, E_l$ be all finite edges of $L$. Let $a_i$ be the lattice length of $E_i$. Then all tropical curves with Newton complex $\nnn$ are parametrized by $a_1, \ldots, a_l> 0$ and $b_1, b_2 \in \rr$. Let $\xxTN( \nnn, \rr^2) \subset \rr^{l+2}$ be the parameter space.

\begin {Prop} \label {Prop:C1}
$\xxTN( \nnn, \rr^2)$ is connected.
\end {Prop}

\begin {proof}
Let $\Gamma_1, \ldots, \Gamma_g$ be all convex cycles of $L$. Let $E_{i(j,1)}, \ldots, E_{i(j,s_j)}$ be all edges of $\Gamma_j$. Let $u_{j,k}$ be the primitive vector of $E_{i(j,k)}$ of positive direction. Then the equation
\begin {equation}
a_{i(j,1)} u_{j,1}+ \cdots + a_{i(j,s_j)} u_{j,s_j}= 0 \label {eq:C1-1}
\end {equation}
is satisfied for any $L \in \xxTN( \nnn, \rr^2)$. Let $H_j \subset \rr^{l+2}$ be the linear subspace defined by equation (\ref{eq:C1-1}). Then
\[
\xxTN( \nnn, \rr^2)= \{ (a_1, \ldots, a_l, b_1, b_2) | a_1, \ldots, a_l> 0 \} \cap (H_1 \cap \cdots \cap H_g).
\]
Thus $\xxTN( \nnn, \rr^2)$ is a relatively open convex cone in $\rr^{l+2}$, which is connected.
\end {proof}

For tropical curves $L, L'$, the notation $\xxNewt(L') \subset \xxNewt(L)$ means that any edge of $\xxNewt( L')$ is an edge of $\xxNewt(L)$.

\begin {Def}
A tropical curve $L'$ is a {\em degeneration} of $L$ if $\Delta(L')= \Delta(L)$ and $\xxNewt(L') \subset \xxNewt(L)$.
\end {Def}

The set of all degenerations of $L$ is parametrized by $\overline{\xxTN( \nnn, \rr^2)}$. If a Newton polygon $\Delta$ is fixed, all tropical curves have a common degeneration, which is a tropical curve consisting of one vertex. The set of the Newton polygons is countably infinite. Therefore, the space $\xxTN( \rr^2)$ of all tropical curves is a disjoint union of countably many closed cones in affine spaces.

\begin {Cor} \label {Prop:D1}
Tropical curves $L, L' \in \xxTN( \rr^2)$ lie in the same connected component if and only if $\Delta(L)= \Delta(L')$.
\end {Cor}

\section {Proof of Theorem \ref{Prop:A1}} \label {Sec:E1}

Let $\Lambda$ be the cycle in $C$. ($C$ is a graph with first Betti number $b_1( C)= 1$. So there is a unique cycle.) $\Lambda$ is equipped with the positive direction defined as follows. For any edge $E \subset \Lambda$ with vertices $V_1, V_2$, the direction from $V_1$ to $V_2$ is positive if and only if
\[
\dvec{V_1 V_2} \times u_{E'} < 0,
\]
where $E'$ is any edge of any tropical curve $L$ that intersects $E$ transversely, and $u_{E'}$ is the weighted primitive vector of $E'$ starting at the vertex inside $\Lambda$.
\par
Let $\oo$ be any point of $\Lambda$. We consider the map $\pi \colon \rr \rightarrow \Lambda$ with the following properties. \\
i) $\pi(0)= \oo$. \\
ii) $\pi$ is increasing with respect to the positive direction of $\Lambda$. \\
iii) $\pi$ is compatible with the lattice length. i.e. If $\pi [a, b]$ ($a, b \in \rr$) is contained in an edge of primitive vector $u$, then
\[
\xxlength( \pi [a, b])=| u|( b- a).
\]
\par
From Lemma \ref{Prop:A5} and Lemma \ref{Prop:A6}, $\xxT(C)$ is generated by $\{ P- \oo| P \in \Lambda \}$. $\Lambda$ has the group structure induced by $\pi$. The following lemma \ref{Prop:A7} says that the map
\[
\ph \colon \Lambda \rmap \xxT(C)
\]
is a surjective homomorphism. \par

\begin {proof}[Proof of that $\ph$ is a homomorphism]
Let $a, b \in \rr$. We have $\pi(a+ b)- \pi(b)= \pi(a)- \oo$ from lemma \ref{Prop:A7}. $\pi(a+ b)- \oo= ( \pi(a)- \oo) +( \pi(b)- \oo)$. Thus $\ph(a+ b)= \ph(a)+ \ph(b)$.
\end {proof}

\begin {Lem} \label {Prop:A7}
If $\pi ( a)= P$, $\pi ( a')= P'$, $\pi ( b)= Q$, $\pi ( b')= Q'$, $a'- a= b'- b$, then $P'- P \sim Q'- Q$.
\end {Lem}

\begin {proof}
(Figure \ref{Fig:A7}) We may assume that $a'- a> 0$ is small enough so that $P, P'$ lie on the same edge $E$, and $Q, Q'$ on $F$. We may also assume that $E, F$ are adjacent at a common vertex $R$, and $E \not= F$. From Lemma \ref{Prop:A3}, we may assume that $P, P', Q, Q'$ are interior points of edges, and $[P, Q']$ has rational slope. \par
Let $L$ be the line passing through $P, Q'$. Then
\[
E \cdot L= |u_E \times u_L| P, \]\[
F \cdot L= |u_F \times u_L| Q.
\]
Let $P'', Q''$ be points of $E, F$ such that
\[
\ddvec{P P''}= \frac{1}{|u_E \times u_L|} \ddvec{P P'}, \]\[
\ddvec{Q' Q''}= \frac{1}{|u_F \times u_L|} \ddvec{Q' Q}.
\]
Then
\begin {eqnarray*}
|\ddvec{P P''} \times u_L| & = & \frac{| \ddvec{P P''}|}{| u_E|} | u_E \times u_L| \\
& = & \frac{| \ddvec{P P'}|}{| u_E|} \\
& = & \frac{| \ddvec{Q' Q}|}{| u_F|} \quad ( \mbox{because } a'- a= b'- b) \\
& = & |\ddvec{Q' Q''} \times u_L|,
\end {eqnarray*}
which means that $[P'', Q'']$ is parallel to $L$. Thus
\[
| u_E \times u_L|( P''- P) \sim |u_F \times u_L|( Q'- Q'').
\]
The interval $[P, P']$ is divided into $| u_E \times u_L|$-number of intervals with the same length. From lemma \ref{Prop:A3}, we have $P'- P \sim | u_E \times u_L|( P''- P)$. Thus $P'- P \sim Q'- Q$.
\end {proof}

\begin {figure}[ht]
\begin {center} \input {./picture/bhA7.tex} \end {center}
\caption {}
\label {Fig:A7}
\end {figure}

Let $E_1, \ldots, E_N$ be all edges of $\Lambda$ ordered by the positive direction. Let $\lambda \colon C \rightarrow \Lambda$ be the canonical surjection. For a tropical curve $L \in \xxTN( \rr^2)$, we define $\sigma(L) \in \Lambda$ as follows.
\[
C \cdot L= P_1+ \cdots+ P_r, \]\[
\sigma(L)= \lambda(P_1)+ \cdots + \lambda(P_r).
\]

\begin {Lem} \label {Prop:C2}
$\sigma \colon \xxTN( \rr^2) \rightarrow \Lambda$ is locally constant.
\end {Lem}

\begin {proof}
$\sigma$ is continuous by definition of the stable intersection. Let $\{ L_t | 0 \leq t \leq 1 \}$ be a continuous family of tropical curves with Newton complex $\nnn$ such that $C$ intersects $L_t$ transversely for any $t$. Let $P_{ijt} \in E_i$ be the points such that $E_i \cdot L_t= \sum_{j} P_{ijt}$. Then there are edges $L_{ijt} \subset L_t$, and vectors $u_{ij} \in \rr^2$ such that \\
i) $E_i \cap L_{ijt}= P_{ijt}$, \\
ii) $u_{ij}$ is the weighted primitive vector of $L_{ijt}$ starting at the vertex inside $\Lambda$. \\
Let $\mu_{ij}$ be the multiplicity $\mu_{P_{ijt}}$ (which is constant for $t$). For $L_0$ and $L_1$, we have the moment condition inside $\Lambda$:
\[
\sum_{i,j} \xxmoment(u_{ij} / \mu_{ij}, P_{ij0})= 0, \]\[
\sum_{i,j} \xxmoment(u_{ij} / \mu_{ij}, P_{ij1})= 0.
\]
From these,
\[
\sum_{i,j} \left( \ddvec{P_{ij0} P _{ij1}} \times u_{ij} / \mu_{ij} \right)= 0.
\]
Since $u_{E_i} \times u_{ij}= - \mu_{ij}$, this means
\[
\sum_{i,j} \left( \lambda(P_{ij0})- \lambda(P_{ij1}) \right)= 0.
\]
Thus $\sigma(L_0)= \sigma(L_1)$.
\end {proof}

\begin {proof}[Proof of the injectivity of $\ph$]
Suppose
\[
P- Q= C \cdot L- C \cdot L', \]\[
P, Q \in \Lambda, \]\[
\Delta(L)= \Delta(L').
\]
From Corollary \ref{Prop:D1} and Lemma \ref{Prop:C2}, we have $\sigma(L)= \sigma(L')$. Thus $P= Q$.
\end {proof}

\section {Proof of Theorem \ref{Thm:E1}} \label{Sec:F1}

Let $H$ be the subgroup of $\xxDiv^0(C)$ generated by all $(P'- P)- (Q'- Q)$ that satisfies the assumption of Lemma \ref{Prop:A3}. We have a homomorphism
\[
\xxDiv^0(C)/ H \rmap \xxT(C).
\]
\par
Let $E_1, \ldots, E_N$ be all finite edges of $C$. Let $V_i, V'_i$ be the vertices of $E_i$. We have a homomorphism
\begin {eqnarray*}
\pi_i \colon \rr & \rmap & \xxDiv^0(C)/ H \\
c & \mapsto & m( P- V_i),
\end {eqnarray*}
where $m> 0$ is any integer such that
\[
\frac{ c}{ m} \leq \xxlength( [V_i, V'_i]),
\]
and $P \in E_i$ is the point such that
\[
\frac{ c}{ m}= \xxlength( [V_i, P]).
\]
This is well-defined by definition of $H$. So we have a homomorphism
\[
\pi \colon \rr^N \rmap \xxDiv^0(C)/ H.
\]
Let $F= \xxKer \pi$. The composite map
\[
\rr^N/ F \rmap \xxDiv^0(C)/ H \rmap \xxT(C)
\]
is surjective because of Lemma \ref{Prop:A5}.
\par
We fix a vertex $\oo \in C$. Let $\lambda \colon C \rightarrow \rr^N/ F$ be the map such that \\
i) $\pi( \lambda(P))= P- \oo$ for $P \in E_i$, \\
ii) every ray of $C$ is contracted by $\lambda$. \\
Let $\widetilde{ \lambda}_i \colon E_i \rightarrow \rr^N$ be a lift of $\lambda |_{E_i} \colon E_i \rightarrow \rr^N/ F$ with the following condition.
\[
\widetilde{ \lambda}_i( P)= \widetilde{ \lambda}_i( V_i)+ \delta e_i,
\]
where
\[
\delta= \xxlength( [V_i, P]). 
\]
\par

\begin {Lem} \label {Prop:G1}
$F$ is a free abelian group of rank $g$.
\end {Lem}

\begin {proof}
Let $\Lambda_1, \ldots, \Lambda_g$ be cycles generating $H_1(C, \zz)$. Let $V_{j1}, \ldots, V_{js_j}$ be all vertices of $\Lambda_j$, and take $V_{j0}= V_{js_j}$. Let $a_{jk} \in \rr^N$ be the vector such that $\pi( a_{jk})= V_{jk}- V_{j, k-1}$. Let $F'$ be the abelian group generated by $a_1, \ldots, a_g$, where $a_j= \sum_{k} a_{jk}$. Then we have $F' \subset F$. The canonical map
\[
\lambda' \colon C \rmap \rr^N/ F'
\]
induces a homomorphism
\[
\lambda'' \colon \xxDiv^0(C)/ H \rightarrow \rr^N/ F'.
\]
So we have $F'= F$.
\end {proof}

For a tropical curve $L \in \xxTN( \rr^2)$, we define $\sigma(L) \in \rr^N/ F$ as follows.
\[
C \cdot L= P_1+ \cdots+ P_r, \]\[
\sigma(L)= \lambda(P_1)+ \cdots + \lambda(P_r).
\]
The map
\[
\sigma \colon \xxTN( \rr^2) \rmap \rr^N/ F
\]
is continuous by definition of the stable intersection.
\par
Given a continuous path $\gamma \colon [0, 1] \rightarrow \xxTN( \rr^2)$, let
\[
\widetilde{ \sigma}_{ \gamma} \colon [0, 1] \rightarrow \rr^N
\]
be a lift of $\sigma |_{ \gamma} \colon [0, 1] \rightarrow \rr^N/ F$, and let
\[
v( \gamma)= \widetilde{ \sigma}_{ \gamma}( 1)- \widetilde{ \sigma}_{ \gamma}( 0).
\]
Let $W \subset \rr^N$ be the set of all $v( \gamma)$. 

\begin {Lem} \label {Prop:G2}
$W$ is a linear subspace in $\rr^N$.
\end {Lem}

\begin {proof}
$W$ is a subgroup of $\rr^N$. Indeed, $v( \gamma)+ v( \gamma') \in W$, because $\sigma(L \cup L')= \sigma(L)+ \sigma(L')$.
\par
Recall that each component of $\xxTN( \rr^2)$ is a closed cone in an affine space. Given a Newton complex $\nnn$, let
\[
\xxTN( \nnn, C)= \left\{ L \in \xxTN( \rr^2) | \xxNewt(C \cup L)= \nnn \right\}.
\]
Every element of $W$ is written as $v( \gamma_1)+ \cdots+ v( \gamma_s)$ for linear paths
\[
\gamma_j \colon [0, 1] \rightarrow \overline{ \xxTN( \nnn_j, C)}
\]
and Newton complexes $\nnn_1, \ldots, \nnn_s$. For $0 \leq t \leq 1$, we have
\[
t \left( v( \gamma_1)+ \cdots+ v( \gamma_s) \right)= v( \gamma_{1, t})+ \cdots+ v( \gamma_{s, t}),
\]
where $\gamma_{j, t} \colon [0, t] \rightarrow \overline{ \xxTN( \nnn_j, C)}$ is a restriction of $\gamma_j$. So $W$ is a linear subspace in $\rr^N$.
\end {proof}

By Corollary \ref{Prop:D1}, $W$ is the kernel of the surjection $\rr^N/ F \rightarrow \xxT(C)$. So we have an isomorphism
\[
\rr^N/( F+ W) \cong \xxT(C).
\]
\par
Note that $\rr^N/ W$ is a vector space generated by $e_1, \ldots, e_N$. Theorem \ref{Thm:E1} follows from the following lemma.

\begin {Lem} \label {Prop:E1}
Suppose that $\ccx := \overline{C \setminus( E_1 \cup \cdots \cup E_g)}$ is a maximal tree in $C$. Then $\langle e_1, \ldots, e_g \rangle$ is a basis of $\rr^N/ W$ over $\rr$.
\end {Lem}

\begin {proof}[Proof of that $\langle e_1, \ldots, e_g \rangle$ is a generator]
For $g+1 \leq j \leq N$, there is a piecewise linear curve $M \subset \rr^2$ containing two rays such that $\ccx \cap M$ is a single point of $E_j$. The following lemma says that $e_j \in \qq \langle e_1, \ldots, e_g \rangle$ in $\rr^N/ W$.
\end {proof}

In the proof of the following lemma, the notation $[P, Q, \infty)$ means the ray with vertex $P$ that contains $Q$.

\begin {Lem} \label {Prop:F1}
Let $S, T$ be rays in $\rr^2$ such that
\[
S \cap C= \emptyset, \quad T \cap C= \emptyset.
\]
Let $P_0, \ldots, P_r, P'_0, \ldots, P'_r \in \rr^2$ be rational points such that
\[
P_0, P'_0 \in S, \quad P_r, P'_r \in T, \]\[
[P_i, P'_i] \cap C= \emptyset, \]\[
\ddvec{P_{i-1} P_i} \mbox{ and } \ddvec{P'_{i-1} P'_i} \mbox{ have the same direction}.
\]
Then there are integers $m_1, \ldots, m_r> 0$ such that
\[
\sum_{i=1}^{r} m_i \left( C \cdot [P_{i-1}, P_i] \right) \sim \sum_{i=1}^{r} m_i \left( C \cdot [P'_{i-1}, P'_i] \right).
\]
\end {Lem}

\begin {proof}
Let $L$ be the tropical curve with vertex $P_i$, consisting of three rays $L_0, L_1, L_2$ such that
\[
P_{i-1} \in L_1, \quad P_{i+1} \in L_2, \]\[
[P_i, P'_i] \mbox{ is parallel to } L_0.
\]
Let $w_i= \xxwt(L_1)$, $v_i= \xxwt(L_2)$. Move $L$ by $\dvec{P_i P'_i}$, and denote it by $L'$. Comparing $C \cdot L$ and $C \cdot L'$, we have
\[
C \cdot \left( w_i [P_i, P_{i-1}, \infty) \cup v_i [P_i, P_{i+1}, \infty) \right) \sim C \cdot \left( w_i [P'_i, P'_{i-1}, \infty) \cup v_i [P'_i, P'_{i+1}, \infty) \right)
\]
for $1 \leq i \leq r-1$. Similarly,
\[
C \cdot \left( w_0 [P_0, S) \cup v_0 [P_0, P_1, \infty) \right) \sim C \cdot \left( w_0 [P'_0, S) \cup v_0 [P'_0, P'_1, \infty) \right), \]\[
C \cdot \left( w_r [P_r, P_{r-1}, \infty) \cup v_r [P_r, T) \right) \sim C \cdot \left( w_r [P'_r, P'_{r-1}, \infty) \cup v_r [P'_r, T) \right).
\]
Taking a suitable linear combination and canceling lines, we have the statement.
\end {proof}

\begin {proof}[Proof of that $\langle e_1, \ldots, e_g \rangle$ are linearly independent]
Suppose
\[
a_1 e_1+ \cdots+ a_g e_g= v( \gamma), \]\[
\gamma \colon [0, 1] \rmap \xxTN( \rr^2).
\]
For $1 \leq j \leq g$, there is a piecewise linear closed curve $\Lambda$ in $\ccx \cup E_j$. Let $E_{i(0)}, \ldots, E_{i(s)}$ be all edges of $\Lambda$, and take $E_{i(0)}= E_j$. Let $\e_k \in \{1, -1 \}$ be the sign of $E_{i(k)}$ with respect to the positive direction of $\Lambda$: The sign is positive if the coordinate of $\rr^N$ is increasing. Let $h \colon \rr^N \rightarrow \rr$ be the linear map defined by
\[
h= \e_0 \check{e}_{i(0)}+ \cdots+ \e_s \check{e}_{i(s)},
\]
where $\check{e}_i$ denotes the dual basis. Then we have
\[
h( \widetilde{ \sigma}_{ \gamma}( 1))- h( \widetilde{ \sigma}_{ \gamma}( 0))= h( v( \gamma))= \e_0 a_j.
\]
To complete the proof, we show that the map $h \circ \widetilde{ \sigma}_{ \gamma} \colon [0, 1] \rightarrow \rr$ is constant for any continuous path $\gamma$. \par
We may assume that $C$ intersects $L_t := \gamma(t)$ transversely for any $t$. There are points $P_{kmt} \in E_{i(k)}$, edges $L_{kmt} \subset L_t$, and vectors $u_{km} \in \rr^2$ such that \\
i) $E_{i(k)} \cdot L_t= \sum_{m} P_{kmt}$, \\
ii) $E_{i(k)} \cap L_{kmt}= P_{kmt}$, \\
iii) $u_{km}$ is the weighted primitive vector of $L_{kmt}$ starting at the vertex inside $\Lambda$. \\
Let $\mu_{km}$ be the multiplicity $\mu_{P_{kmt}}$ (which is constant for $t$). For $L_0$ and $L_1$, we have the moment condition inside $\Lambda$:
\[
\sum_{k,m} \xxmoment(u_{km} / \mu_{km}, P_{km0})= 0, \]\[
\sum_{k,m} \xxmoment(u_{km} / \mu_{km}, P_{km1})= 0.
\]
From these,
\[
\sum_{k,m} \left( \ddvec{P_{km0} P _{km1}} \times u_{km} / \mu_{km} \right)= 0.
\]
Let $u_k$ be the primitive vector of $E_{i(k)}$ starting at $V_{i(k)}$. Since $u_k \times u_{km}=- \e_k \mu_{km}$, we have
\[
\sum_{k,m} h \left( \widetilde{ \lambda}_{i(k)}(P_{km0})- \widetilde{ \lambda}_{i(k)}(P_{km1}) \right)= 0.
\]
Thus $h( \widetilde{ \sigma}_{ \gamma}( 0))= h( \widetilde{ \sigma}_{ \gamma}( 1))$.
\end {proof}

\end {document}

%% file: picture/bhA3_b.tex
\begin {picture}(0, 0)
\includegraphics {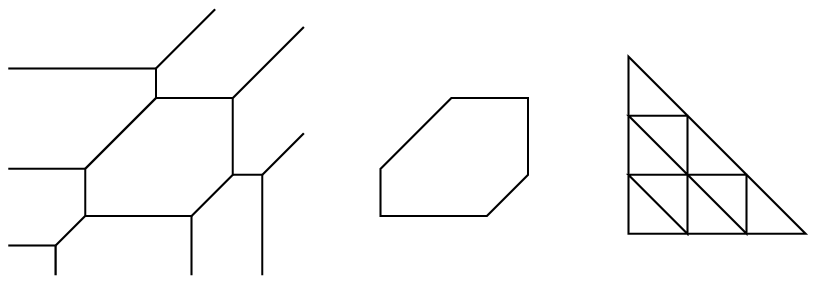}
\end {picture}
\begin {picture}(350, 150)
\unitlength= 1 mm

\put (25, -5) {$C$}
\put (67, -5) {$\Lambda$}
\put (111, -5) {$\xxNewt(C)$}

\end {picture}

%% file: picture/bhA2.tex
\begin {picture}(0, 0)
\includegraphics {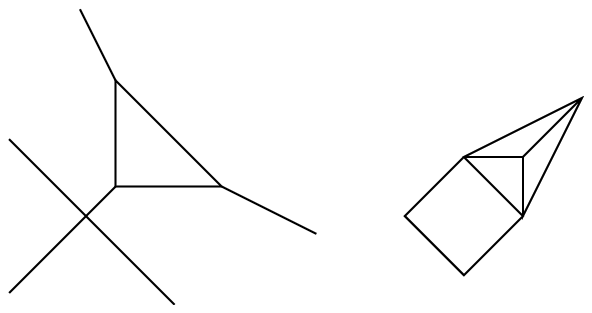}
\end {picture}
\begin {picture}(350, 150)
\unitlength= 1 mm

\put (28, -5) {$C$}
\put (82, -5) {$\xxNewt(C)$}

\end {picture}

%% file: picture/bhB1.tex
\begin {picture}(0, 0)
\includegraphics {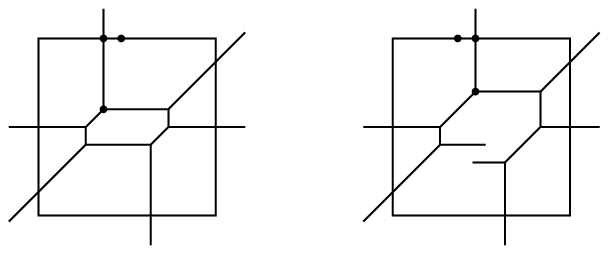}
\end {picture}
\begin {picture}(350, 150)
\unitlength= 1 mm

\put (10, 36) {$P$}
\put (70, 36) {$P$}
\put (19, 36) {$P'$}
\put (79, 36) {$P'$}
\put (13, 42) {$E_1$}
\put (76, 42) {$E'_1$}
\put (7, 24) {$V_{E_1}$}
\put (70, 27) {$V_{E'_1}$}

\end {picture}

%% file: picture/bhA5.tex
\begin {picture}(0, 0)
\includegraphics {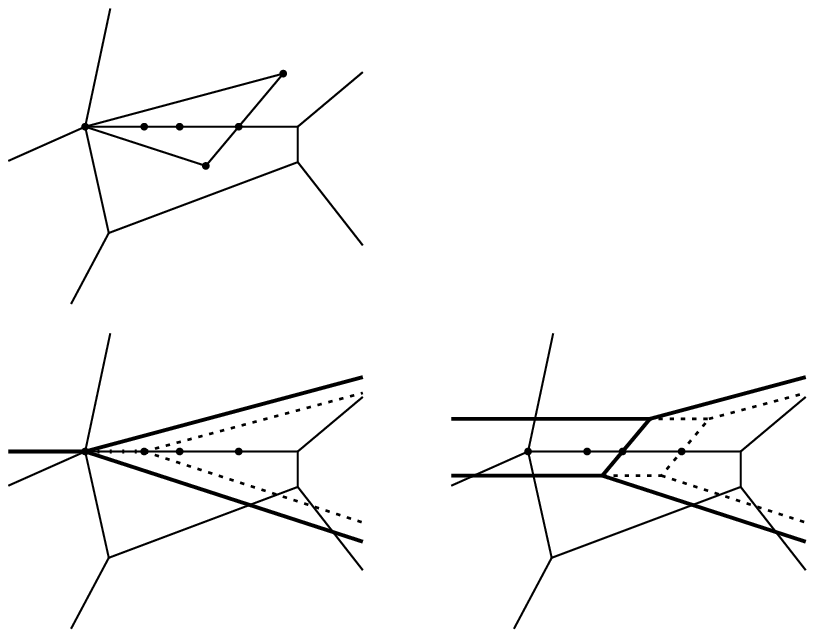}
\end {picture}
\begin {picture}(350, 300)
\unitlength= 1 mm

\put (9, 87) {$P$}
\put (18, 89) {$P'$}
\put (32, 87) {$Q$}
\put (42, 87) {$Q'$}
\put (27, 76) {$R_1$}
\put (44, 97) {$R_2$}

\put (9, 32) {$P$}
\put (18, 34) {$P'$}
\put (1, 32) {$L$}

\put (107, 32) {$Q$}
\put (117, 32) {$Q'$}
\put (76, 37) {$M$}

\end {picture}

%% file: picture/bhA6.tex
\begin {picture}(0, 0)
\includegraphics {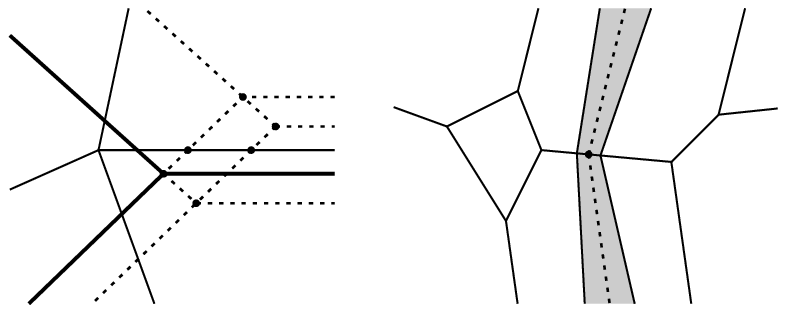}
\end {picture}
\begin {picture}(350, 150)
\unitlength= 1 mm

\put (44, 32) {$R_1$}
\put (39, 37) {$R_2$}
\put (31, 13) {$R_3$}
\put (23, 17) {$R_4$}
\put (3, 44) {$M_4$}
\put (55, 26) {$E$}
\put (39, 23) {$P$}
\put (25, 29) {$Q$}

\put (100, -4) {$\theta_2$}
\put (102, 53) {$\theta_1$}
\put (101, 27) {$P$}

\end {picture}

%% file: picture/bhA7.tex
\begin {picture}(0, 0)
\includegraphics {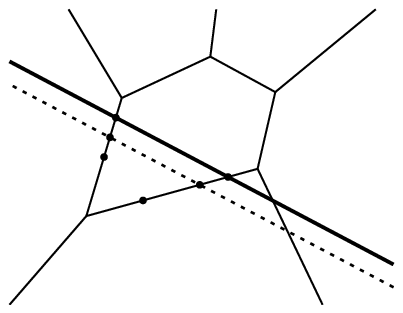}
\end {picture}
\begin {picture}(200, 150)
\unitlength= 1 mm

\put (65, 7) {$L$}
\put (19, 31) {$P$}
\put (7, 26) {$P''$}
\put (9, 21) {$P'$}

\put (19, 12) {$Q$}
\put (28, 15) {$Q''$}
\put (34, 24) {$Q'$}

\end {picture}

%% file: TJACk6.bbl
\begin {thebibliography}{9}

\bibitem {Gath}
A. Gathmann.
\newblock Tropical algebraic geometry.
\newblock Preprint, arXiv:math.AG/0601322.

\bibitem {Mik3}
Grigory Mikhalkin.
\newblock Decomposition into pairs-of-pants for complex algebraic hypersurfaces.
\newblock Topology 43, Issue 5 (2004), 1035-1065.

\bibitem {Mik}
Grigory Mikhalkin.
\newblock Enumerative tropical algebraic geometry in $\rr^2$.
\newblock J. Amer. Math. Soc. 18 (2005), 313-377.

\bibitem {Mik2}
G. Mikhalkin and I. Zharkov.
\newblock Tropical curves, their Jacobians and Theta functions.
\newblock Preprint, arXiv:math.AG/0612267.

\bibitem {Stu2002}
Bernd Sturmfels.
\newblock Solving systems of polynomial equations.
\newblock Volume 97 of CBMS Regional Conference Series in Mathematics (2002).

\bibitem {First Steps}
J. Richter-Gebert, B. Sturmfels, and T. Theobald.
\newblock First Steps in Tropical Geometry.
\newblock Preprint, arXiv:math.AG/0306366.

\bibitem {elliptic curve}
M. D. Vigeland.
\newblock The group law on a tropical elliptic curve.
\newblock Mathematica scandinavica 104, issue 2 (2009), 188-204.




\end {thebibliography}
